\documentclass[12pt,draft]{article}   

\usepackage{amsthm,amssymb,amsmath,amscd,amstext,mathrsfs,latexsym}

\newtheorem{definition}{Definition}[section]
\newtheorem{lemma}{Lemma}[section]

\newtheorem*{remark}{{\it Remark}}

\newcommand{\nc}{\newcommand}
\nc{\N}{{\mathbb N}}
\nc{\Z}{{\mathbb Z}}
\nc{\R}{{\mathbb R}}
\nc{\C}{{\mathbb C}}
\nc{\HH}{{\mathbb H}}
\nc{\ca}{{\mathscr A}}
\nc{\cb}{{\mathscr B}}
\nc{\cu}{{\mathscr U}(1)}
\nc{\ch}{{\mathscr H}}
\nc{\dd}{{\rm d}}
\nc{\DD}{{\rm D}}
\nc{\ii}{{\bf i}}

\begin{document}

\title{$S$-duality in Abelian gauge theory revisited}

\author{G\'abor Etesi\\
\small{{\it Department of Geometry, Mathematical Institute, Faculty of
Science,}}\\
\small{{\it Budapest University of Technology and Economics,}}\\
\small{{\it Egry J. u. 1, H \'ep., H-1111 Budapest, Hungary
\footnote{e-mail: {\tt etesi@math.bme.hu}}}}\\
\'Akos Nagy\\
\small{{\it Budapest University of Technology and Economics}}
\footnote{e-mail: {\tt nagyak@math.bme.hu}}}

\maketitle

\pagestyle{myheadings}
\markright{G. Etesi, \'A. Nagy: Abelian $S$-duality revisited}

\thispagestyle{empty}

\begin{abstract}
Definition of the partition function of ${\rm U}(1)$ gauge theory is 
extended to a class of four-manifolds containing all compact spaces and 
the asymptotically locally flat (ALF) ones including the 
multi-Taub--NUT sopaces. The partition function is 
calculated via zeta-function regularization and heat kernel 
techniques with special attention to its modular properties. 

In the compact case, compared with the purely topological result of 
Witten, we find a non-trivial curvature correction to the modular weights of 
the partition function. But $S$-duality can be restored by adding 
gravitational counter terms to the Lagrangian in the usual way. 

In the ALF case however we encounter non-trivial difficulties stemming 
from original non-compact ALF phenomena. Fortunately our careful definition 
of the partition function makes it possible to circumnavigate them and 
conclude that the partition function has the same modular properties as in 
the compact case.
\end{abstract}

\centerline{AMS Classification: Primary: 81T13; Secondary: 81Q30, 
57M50, 11F37, 35K08}
\centerline{Keywords: {\it S-duality; Partition 
function; $L^2$ cohomology; Zeta-function regularization; Heat kernel}}


\section{Introduction}
\label{one}


The long standing conjecture asserts that quantum gauge theory has a 
symmetry exchanging strong and weak coupling as well as electric and magnetic 
fields. This conjecture originated with the work of Montonen and Olive 
\cite{mon-oli} from 1977 who proposed a symmetry in quantum gauge theory with 
the above properties and also interchanging the gauge group $G$ with its dual 
group $G^\vee$. It was soon realized however that this duality is more 
likely to hold in an $N=4$ supersymmetrized theory \cite{oli-wit}. 

The original Montonen--Olive conjecture proposed a $\Z_2$ symmetry 
exchanging the electric and magnetic charges however in $N=4$ theory it 
is naturally extended to an ${\rm SL}(2,\Z )$ symmetry acting on the complex 
coupling constant
\begin{equation}
\tau :=\frac{\theta}{2\pi}+\frac{4\pi}{e^2}\ii\:\:\in\C^+
\label{tau}
\end{equation}
combining the gauge coupling $e$ and the $\theta$ parameter of the $N=4$ 
theory. In this framework the conjecture can be formulated as follows 
\cite{vaf-wit}.

We say that a not necessarily holomorphic function 
$f:\C^+\rightarrow\C$ on the upper half-plane is an {\it unrestricted modular 
form of weight $(\alpha ,\beta )$} (conventionally supposed to be integers) 
if with respect to 
$\begin{pmatrix}a & b\\ c & d\end{pmatrix}\in{\rm SL}(2,\Z )$ it transforms as
\[f\left(\frac{a\tau +b}{c\tau +d}\right) =(c\tau +d)^\alpha 
(c\overline{\tau} +d)^\beta f(\tau ).\]
If $\alpha\not= 0$ or $\beta\not= 0$ we also say sometimes that a {\it 
modular anomaly} is present in $f$. Since ${\rm SL}(2,\Z )$ is generated 
by $T=\begin{pmatrix}1 & 1\\ 0 & 1\end{pmatrix}$ and 
$S=\begin{pmatrix}0 & -1\\1 & 0\end{pmatrix}$ 
modular properties are sufficient to be checked under the simple 
transformations 
$f(\tau )\mapsto f(\tau +1)$ and $f(\tau )\mapsto f(-1/\tau )$. Modular 
forms play an important role in classical number theory \cite{kna,ser}. 

The {\it electric-magnetic duality conjecture} in its simplest 
form asserts that over a four-manifold $(M,g)$ the partition function of 
a (twisted) $N=4$ supersymmetric quantum gauge theory with 
simply-laced gauge group $G$ is modular in the sense that it satisfies
\[\left\{\begin{array}{ll}
                     \mbox{$Z(M,g,G,\tau +1)=Z(M,g,G,\tau )$
                       \:\:\:\:\:(``level $1$ property'')}\\
                     \mbox{$Z(M,g,G,-1/\tau )=Z(M,g,G^\vee ,\tau )$
                      \:\:\:\:\:(``$S$-duality property'')}
\end{array}\right.\]
where $G^\vee$ is the Goddard--Nuyts--Olive or Langlands dual group to 
$G$. The first symmetry is a classical one while the second is expected to 
reflect the true quantum nature of gauge theories in the 
sense that it connects two theories on the quantum 
level which are classically different. In general $G$ is not isomorphic 
to its dual however for example if $G\cong{\rm U}(1)$ then it is. 
Moreover if $\theta =0$ this duality reduces to the 
original Montonen--Olive conjecture. 

We do not attempt here to survey the long, diverse and colourful history 
of the conjecture and its variants rather refer to the introductions of 
\cite{vaf-wit, kap-wit}. We just mention that for instance it led to 
highly non-trivial predictions about the number of $L^2$ harmonic forms 
on complete manifolds due to Sen \cite{sen} (cf. also \cite{gib, 
hau-hun-maz,seg-sel}) and the latest chapter of 
the story relates the electric-magnetic duality conjecture with the 
geometric Langlands program of algebraic geometry due to Kapustin 
and Witten \cite{kap-wit}, see also \cite{fre}.

In this work, motivated by papers of Witten \cite{wit1, wit2} we focus 
attention to the modular properties of the partition function 
\[Z(M,g,\tau ):=Z(M,g,{\rm U}(1),\tau )\]
of (supersymmetric) Abelian or ${\rm U}(1)$ or Maxwell gauge theory over 
various four-manifolds because in this case the calculations can be carried 
out almost rigorously. 

The paper is organized as follows. In Sect. \ref{two} we offer an extended 
definition of the partition function of Abelian gauge theory (cf. Eq. 
(\ref{particiofv}) here) such that it gives back Witten's 
\cite{wit1} in case of compact four-manifolds moreover the 
definition continues to make sense for a certain class of non-compact 
geometries, the so-called {\it asymptotically locally flat} 
(ALF) spaces including the flat $\R^3\times S^1$, the multi-Taub--NUT 
spaces and the Riemannian Schwarzschild and Kerr solutions, etc. The 
subtle point of this extended definition is imposing a natural boundary 
condition at infinity on finite action classical solutions of Maxwell theory 
in order to save the modular properties of the partition function. This 
boundary condition is the so-called {\it strong holonomy condition} and is 
used to rule out certain classical solutions having non-integer energy. 
It has also appeared already in the approach to SU(2) instanton moduli 
spaces over ALF geometries \cite{ete, ete-jar}.

Then in Sect. \ref{three} we repeat Witten's calculation \cite{wit1} of  
$Z(M,g,\tau )$ in the compact case. After summing up over a discrete 
set leading to $\vartheta$-functions the partition function is given by a 
formal infinite dimensional integral and can be calculated quite rigorously 
via $\zeta$-function regularization. On the way we formulate a natural 
``Fubini principle'' for these formal integrals which converts 
them into successive integrals. We come up with an expression involving 
the zero values of various $\zeta$-functions and their derivatives. 
The zero values of these $\zeta$-functions can be calculated by the aid of 
heat kernel techniques. With the help of these tools we find that 
\begin{equation}
\left\{\begin{array}{ll}
Z(M,g,\tau +2)= Z(M,g,\tau )\\
Z(M,g,-1/\tau )=(-\ii)^{\frac{1}{2}\sigma (M)}
\:\tau^\alpha\:\overline{\tau}^\beta \:Z(M,g,\tau )\end{array}\right.
\label{eredmeny}
\end{equation}
i.e., up to a factor it is a level $2$ modular form but with 
non-integer weights (cf. Eq. (\ref{sulyok}) here)
\[\alpha =\frac{1}{4}\left(\chi (M)+\sigma (M)+\mbox{curvature 
corrections}\right)\] 
and 
\[\beta =\frac{1}{4}\left(\chi (M)-\sigma (M)+\mbox{curvature 
corrections}\right)\]
yielding that the modular weights are not purely topological as claimed 
in \cite{wit1, wit2}. As a consequence $S$-duality fails in its 
simplest form but it can be saved by adding usual 
gravitational $c$-number terms to the naive Lagrangian of Maxwell theory 
on a curved background \cite{wit1}.

Before proceeding to the non-compact calculation we make a digression 
in Sect. \ref{four} and clarify the role played by the strong holonomy 
condition. We will see that without it (i.e., simply taking the 
definition of the partition function from the compact case) the partition 
function for instance over the multi-Taub--NUT spaces would seriously fail to 
be modular as a consequence of the presence of classical 
solutions with continuous energy spectrum.

Finally as a novelty in Sect. \ref{five} we repeat the calculation over 
ALF spaces. $S$-duality over non-compact geometries is less known 
(cf. the case of ALE spaces in \cite{vaf-wit}). However before obtaining 
some results we have to overcome another 
technical difficulty caused by non-compactness. Namely if one wishes to 
calculate the formal integrals by $\zeta$-function regularization again then 
first one has to face the fact that the spectra of various 
differential operators are continuous hence the existence of their 
$\zeta$-functions is not straightforward. Consequently we take a 
truncation of the original manifold and impose Dirichlet boundary condition 
on the boundary. This boundary condition is compatible with the 
finite action assumption on connections. After defining everything 
correctly in this framework, we let the boundary go to infinity. We will find 
that the modular weights converge in this limit and give back a formula very 
similar to the compact case above. The only difference is that the various 
Betti numbers which enter the modular weights are mixtures of true $L^2$ Betti 
numbers and limits of Dirichlet--Betti numbers which are remnants of the 
boundary condition (see Eq. (\ref{alfsulyok}) here). 

Consequently the partition function transforms akin to (\ref{eredmeny}) 
again however its modular anomaly cannot be cancelled by adding $c$-number 
terms. However up to a mild topological condition on the infinity of an ALF 
space (cf. Eq. (\ref{perem}) here) all these Betti numbers can be converted 
into $L^2$ ones again and the cancellation of the modular anomaly goes as in 
the compact case hence $S$-duality can be saved. It is interesting that while 
the Riemannian Schwarzschild and Kerr not, the multi-Taub--NUT family 
satisfies this topological condition. 

We clarify at this point that throughout this paper questions related to 
the contributions of various determinants to the partition function will 
be suppressed.

We close this introduction by making a comment about the way of 
separating rigorous mathematical steps from intuitive ones in the text. To 
help the reader we decided to formulate every rigorous steps in the form of a 
lemma (with proof) or with a clear reference to the literature while the 
other considerations just appear continuously in the text. 
\vspace{0.1in}

\noindent{\bf Acknowledgement.} The authors are grateful to Gyula Lakos 
and Szil\'ard Szab\'o for the stimulating discussions. The first author 
was partially supported by OTKA grant No. NK81203 (Hungary).
 

\section{Preliminary calculations}
\label{two}


In this section we formulate our problem as precisely as possible. Let 
$(M,g)$ be a connected, oriented, complete Riemannian four-manifold without 
boundary. $M$ can be either compact or non-compact. In the non-compact 
case we require $(M,g)$ to have infinite volume. 

Let $L$ be a smooth complex line bundle over $M$. The isomorphism 
classes of these bundles are classified by the elements of the 
group $H^2(M;\Z )$ via their first Chern class $c_1(L)\in H^2(M;\Z )$. 
Putting a Hermitian structure onto $L$ we can pick a 
${\rm U}(1)$-connection $\nabla$ with curvature $F_\nabla$. 
Under ${\mathfrak u}(1)\cong\ii\R$ the curvature is an $\ii\R$-valued 
$2$-form closed by Bianchi identity hence $[\frac{F_\nabla}{2\pi\ii}]\in 
H^2(M;\R )$ is a de Rham cohomology class taking its values in the integer 
lattice. This cohomology class (the Chern--Weil class of $L$) 
is independent of the connection but it is a slightly weaker invariant of 
$L$ namely it characterizes it up to a flat line bundle only. Note that flat 
line bundles over $M$ satisfy $c_1(L)\in{\rm Tor}(H^2(M;\Z ))$.

With real constants $e$ and $\theta$ the usual action of 
(Euclidean) electrodynamics extended with the so-called $\theta$-term 
over $(M,g)$ looks like

\begin{equation}
S(\nabla ,e,\theta ):=-\frac{1}{2e^2}\int\limits_M F_\nabla\wedge *F_\nabla 
+\frac{\ii \theta}{16\pi^2}\int\limits_M F_\nabla\wedge F_\nabla .
\label{hatas}
\end{equation}
The $\theta$-term is a characteristic class hence its variation is 
identically zero consequently the Euler--Lagrange equations of this 
theory are just the usual vacuum Maxwell-equations
\begin{equation}
\dd F_\nabla =0,\:\:\:\:\:\delta F_\nabla =0.
\label{maxwell}
\end{equation}

\begin{remark}\rm If one considers the underlying quantum field theory 
then $\theta$ in (\ref{hatas}), being a non-dynamical variable, remains 
well-defined at the full quantum level meanwhile the coupling constant 
$e$ runs. Therefore, in order to keep its meaning at 
the full quantum level we extend (\ref{hatas}), as usual, to an $N=4$ 
supersymmetric theory \cite{oli-wit,vaf-wit}. This extension yields 
additional 
terms to (\ref{hatas}) however their presence do not influence our 
forthcoming calculations therefore we shall not mention them explicitly in 
this paper.
\end{remark}
 
\noindent Rather we introduce the complex 
coupling constant (\ref{tau}) taking its values on the upper half-plane
and supplementing (\ref{hatas}) re-write the relevant part of the action as
\begin{equation}
S(\nabla ,\tau )=
\frac{\ii\pi}{2}\tau\left(\frac{1}{8\pi^2}\int\limits_M (F_\nabla\wedge
*F_\nabla +F_\nabla\wedge F_\nabla )\right)
+\frac{\ii\pi}{2}(-\overline{\tau})\left(\frac{1}{8\pi^2}\int\limits_M
(F_\nabla\wedge *F_\nabla -F_\nabla\wedge F_\nabla )\right) .
\label{bovites}
\end{equation}
For clarity we note that this expression is exactly the same as (\ref{hatas}).

The orientation and the metric on $M$ is used to form various Sobolev 
spaces. Fix a line bundle $L$ over $M$ and a connection $\nabla^0_L$ on it 
such that $\ii F_{\nabla^0_L}\in L^2(M;\Lambda^2M)$. For any integer 
$l\geqq 2$ then set
\[\ca (\nabla^0_L):=\{ \nabla_L\:\vert\:\mbox{$\nabla_L:=\nabla^0_L+a$ 
with $\ii a\in L^2_l(M;\Lambda^1M)$}\} .\]
This is the $L^2_l$ Sobolev space of ${\rm U}(1)$ connections on $L$ 
relative to $\nabla^0_L$. Notice that this is a vector space (not an 
affine space) and if $L_0\cong M\times\C$ is the trivial line bundle 
and $\nabla^0_{L_0}=\dd$ is the trivial flat connection in the trivial 
gauge on it then $\ca (\nabla^0_{L_0})\cong L^2_l(M;\Lambda^1M)$. Furthermore 
write $\cu_L$ for the $L^2_{l+1}$-completion of the space of gauge 
transformations
\[\{ \gamma -{\rm Id}_L\in C^\infty _0(M;{\rm End}L)\:\vert\:
\Vert\gamma -{\rm Id}_L\Vert_{L^2_{l+1}(M)}<+\infty\:;\:
\mbox{$\gamma\in C^\infty (M;{\rm Aut}L)$ a.e.}\}.\]
Under these assumptions the Sobolev multiplication
theorem ensures us that $\ii F_{\nabla_L}\in L^2(M\:;\:\Lambda^2M)$ that 
is, the curvature $2$-form is always in $L^2$ over $(M,g)$.
The space $\ca (\nabla^0_L)$ is acted upon by $\cu_L$ in the usual 
way; the orbit space $\ca (\nabla^0_L)/\cu_L$ of gauge equivalence 
classes with its quotient topology is denoted by $\cb (\nabla^0_L)$ as usual. 

Next we record an obvious decomposition of the action. 
This decomposition plays a crucial role because it underlies the definition of 
the partition function. Suppose that $\nabla^0_L\in\ca (\nabla^0_L)$ is a 
{\it finite action classical solution on $L$} i.e., it satisfies 
(\ref{maxwell}) and has finite action (\ref{hatas}) or (\ref{bovites}); 
moreover let $\nabla_L\in\ca (\nabla^0_L)$ be another connection on the 
same line bundle $L$. Then it follows from (\ref{bovites}) that
\begin{equation}
S(\nabla_L,\tau )=S(\nabla^0_L+a,\tau )=S(\nabla^0_L,\tau )+\frac{{\rm 
Im}\tau}{8\pi}\Vert\dd a\Vert^2_{L^2(M)}
\label{szetszedes}
\end{equation}
where we write $\Vert\dd a\Vert^2_{L^2(M)}=-\int_M\dd a\wedge *\dd a$ 
since $a$ is pure imaginary. This decomposition is a straightforward 
calculation if $M$ is compact while follows from the $L^2$ control on the 
curvature of $\nabla^0_L$ and the $L^2_1$ control on the perturbation 
$a$ if $M$ is non-compact.

Now we turn to clarify what shall we mean by a ``partition function'' in 
this paper. In our attempt to define it we have to be careful because
we want to include the case of non-compact base spaces as well while
we require that our definition should agree with Witten's \cite{wit1,
wit2} in the compact case. The above decomposition of the action 
indicates that finding the partition function i.e., the integration over 
all finite-action connections should be carried out in two steps: (i) 
summation over a certain set of classical solutions $\nabla^0_L$ (this 
will coincide with summation over all line bundles $L$ in the compact case) 
and then (ii) integration over the orbit space $\cb (\nabla^0_{L_0})$ for the 
trivial flat connection $\nabla^0_{L_0}$ in Coulomb gauge on the trivial 
line bundle $L_0$.

We begin with the description of the orbit space. First we 
have a usual gauge fixing lemma.
\begin{lemma}
Let $(M,g)$ be as above i.e., a connected, oriented complete Riemannian
four-manifold which is either closed (compact without boundary) or open
(non-compact without boundary) with infinite volume. Consider the 
trivial line bundle $L_0\cong M\times\C$ on it.

If $\nabla_{L_0} =\dd +A$ with $A\in\ca (\nabla^0_{L_0})$ is any
${\rm U}(1)$-connection on $L_0$ then there exists an $L^2_{l+1}$ gauge
transformation $\gamma :M\rightarrow{\rm U}(1)$ such that the resulting
$L^2_l$ connection $\nabla '_{L_0}=\dd +A'$ satisfies the Coulomb gauge
condition $\delta A'=0$.
\label{coulomblemma}
\end{lemma}

\noindent {\it Proof.} Take a gauge transformation $\gamma ={\rm e}^{\ii f}$
with a function $f:M\rightarrow\R$. We want to
solve the equation $0=\delta A' =\delta (A+\ii\dd f)$ i.e.,
$\triangle_0f=\ii\delta A$. This equation is solvable if and only if
the real-valued function $\ii\delta A$ is orthogonal in $L^2$ to the
cokernel of the scalar Laplacian $\triangle_0=\delta\dd$ on $(M,g)$. 
Note that $\triangle_0$ is formally self-adjoint hence in fact we need
$\ii\delta A\perp_{L^2}{\rm ker}\triangle_0$. However
${\rm ker}\triangle_0\cong\R$
if $M$ is closed hence the result follows from Stokes' theorem. While
${\rm ker}\triangle_0\cong\{ 0\}$ if $M$ is non-compact but complete and
has infinite volume by a theorem of Yau \cite{yau} hence the result also
follows. $\Diamond$
\vspace{0.1in}

\noindent The time has come to 
introduce $L^2$ cohomology groups on $(M,g)$ in the usual way 
\cite{hau-hun-maz}: we say that a 
real $k$-form $\varphi$ belongs to the {\it $k^{{\it th}}$ reduced $L^2$ 
cohomology group $\overline{H}^k_{L^2}(M)$} if and only if
$\Vert\varphi\Vert_{L^2(M)}<+\infty$ as well as $\dd\varphi =0$ and
$\delta\varphi =0$. Over a complete manifold we can equivalently say 
that $\varphi$ has finite $L^2$ norm and is harmonic i.e., 
$\triangle_k\varphi =0$ where $\triangle_k$ is the Laplacian on $k$-forms. 
If $\ch^k_{L^2}(M)$ denotes the space of $L^2$ harmonic forms then by 
definition $\overline{H}^k_{L^2}(M)\cong\ch^k_{L^2}(M)$ consequently in 
this paper we will write $\ch^k_{L^2}(M)$ for these groups. In general we 
call a non-trivial $L^2$ harmonic form on a complete 
Riemannian manifold {\it non-topological} if either it is exact or not 
cohomologous to a compactly supported differential form. Roughly speaking the 
existence of non-topological $L^2$ harmonic forms are not predictable by 
topological means (cf. \cite{seg-sel}). Note that
in the compact case the $L^2$ groups reduce to ordinary de Rham
cohomology groups $H^k(M;\R )$ and of course there are no non-topological 
$L^2$ harmonic forms. Write $b^k_{L^2}(M):=\dim\ch^k_{L^2}(M)$ 
for the corresponding $L^2$ Betti number (if finite). Of course in 
the compact case $b^k_{L^2}(M)=b^k(M)$ are just the ordinary finite
Betti numbers.
 
Most of these groups admit natural interpretations 
in terms of our theory (\ref{hatas}) as we will see shortly. 

\begin{lemma} Let $M$ be an arbitrary connected manifold and denote by 
$\widehat{M}$ the space of gauge equivalence classes of $L^2$ flat 
${\rm U}(1)$-connections on $M$ (also called the {\em character variety} 
of $M$). 

{\rm (i)} There is an identification 
\[\widehat{M}\cong T^{b_1(M)}\times{\rm Tor}(H_1(M;\Z ))\]
where $T^k$ is a $k$-torus.

{\rm (ii)} If $(M,g)$ is a Riemannian four-manifold as in the previous 
lemma then we have a more informative identification 
\[ \widehat{M}\cong\frac{\ch^1_{L^2}(M)}{2\pi\Lambda^1_{L^2}}
\times{\rm Tor}(H_1(M;\Z ))\]
where $\Lambda^1_{L^2}$ is a co-compact integer lattice in $\ch^1_{L^2}(M)$. 

Moreover keep in mind that $\ch^1_{L^2}(M)\cong H^1(M;\R )$ is 
ordinary de Rham cohomology if $M$ is compact. 
\label{laposlemma}
\end{lemma}

\begin{remark}\rm The lemma implies that if $M$ is a connected, oriented, 
non-compact, complete four-manifold with infinite volume then 
$b^1_{L^2}(M)=b^1(M)$ hence there are no non-topological $L^2$ harmonic 
$1$-forms.
\end{remark}

\noindent{\it Proof.} The gauge equivalence classes of flat 
${\rm U}(1)$-connections are classified by the space
\[\widehat{M}\cong {\rm Hom}(\pi_1(M)\:;\:{\rm U}(1))/{\rm ad\: U}(1)
={\rm Hom}(\pi_1(M)\:;\:{\rm U}(1))\]
since ${\rm U}(1)$ is Abelian. Also commutativity implies that 
the commutator group $\left[\pi_1(M),\pi_1(M)\right]$ 
lies in the kernel of any such homomorphism hence $\pi_1(M)$ can be replaced by 
\[\pi_1(M)/\left[\pi_1(M),\pi_1(M)\right]\cong H_1(M;\Z ).\] 
Write $H_1(M;\Z )=\Lambda_1\times{\rm Tor}(H_1(M;\Z ))$ with 
$\Lambda_1\cong\Z^{b_1(M)}$. Homomorphisms $\rho: \Z\rightarrow {\rm U}(1)$ 
are classified by their trace ${\rm tr}(\rho (1))\in S^1\subset\C$ on the 
generator $1\in\Z$ hence the parameter space is $S^1$. The homomorphisms of the 
finite Abelian group $\rho :{\rm Tor}(H_1(M;\Z ))\rightarrow {\rm U}(1)$ are 
certainly classified by its elements. This yields part (i) of the lemma.

Regarding the second part, suppose that $\nabla^0_{L_0}$ is the trivial 
flat ${\rm U}(1)$-connection on the trivial bundle $L_0\cong M\times\C$ in 
the straightforward gauge i.e., when $\nabla^0_{L_0} =\dd$. It follows 
from (\ref{szetszedes}) that its perturbation $\nabla_{L_0} :=\nabla^0_{L_0}+a$
with $a\in\ca (\nabla^0_{L_0})$ is flat if and only if 
$\dd a=0$. Imposing the Coulomb gauge condition $\delta a=0$ (cf. Lemma 
\ref{coulomblemma}) apparently the gauge-inequivalent perturbations are 
parameterized by the group $\ch^1_{L^2}(M)$. Indeed, if two flat 
connections $\dd +a$ and $\dd +b$ on $L_0$ are in Coulomb gauge then it 
easily follows from the proof of Lemma \ref{coulomblemma} that there is no 
gauge transformation between them of the form $\gamma ={\rm e}^{\ii f}$  with a 
single-valued function $M\rightarrow\R$. 

However if $M$ is not simply-connected then there may still exists a gauge 
transformation such that $b=a+\gamma^{-1}\dd\gamma$
with $\gamma ={\rm e}^{\ii F}$ where $F$ is a multi-valued real function 
on $M$ (or a single-valued function on the universal cover 
$\widetilde{M}$). More precisely, if $\ell :[0,1]\rightarrow M$ is a smooth 
loop with base point $x_0\in M$ then 
$\gamma (\ell (0))=\gamma (\ell (1))=\gamma (x_0)$ implies 
\[F(x_0)=F(\ell (1))=F(\ell (0))+2\pi n_\ell =F(x_0)+2\pi n_\ell\] 
with some integer $n_\ell$. Consider a flat connection $\dd +a$ with 
$a\in\ca (\nabla^0_{L_0})$ and also pick a differentiable loop $\ell$ as 
before. Associated with this loop take a gauge transformation 
$\gamma_\ell ={\rm e}^{\ii F_\ell}$ such that $F_\ell$ is supported in a 
tubular neighbourhood of $\ell$ and $F_\ell (\ell (t))=F_\ell (x_0)+2\pi t$. 
Then obviously $\gamma_\ell\in\cu_{L_0}$ and 
\[a+\gamma^{-1}_\ell\dd\gamma_\ell = a +2\pi\ii\:\dd t\]
moreover the new connection depends only on the homotopy class $[\ell ]$. 
Consider the lattice $\Lambda '$ in $\ch^1_{L^2}(M)$ generated 
by the above translations when $[\ell ]\in \pi_1(M,x_0)$ runs over the 
homotopy classes. It follows that the truely gauge-inequivalent flat 
connections are given by the quotient $\ch^1_{L^2}(M)/\Lambda '$.
But comparing this quotient with part (i) of the lemma we conclude that 
$\Lambda '$ is co-compact hence it must coincide with a full integer 
lattice of $\ch^1_{L^2}(M)$. This gives a connected component 
of $\widehat{M}$ in the second picture.

Note that all the flat $L^2$ perturbations of the trivial flat 
connection are connections on the same trivial bundle $L_0$. The second 
integer cohomology has the straightforward decomposition
\begin{equation}
H^2(M;\Z )=\Lambda^2\times {\rm Tor}(H^2(M;\Z ))
\label{H2}
\end{equation}
into its free part $\Lambda^2\cong\Z^{b^2(M)}$ and its torsion. We 
obtain that generic elements $(a,\alpha )\in\widehat{M}$ 
are flat connections of the shape $\nabla^0_L+a$ where $\nabla^0_L$ is a fixed 
flat connection on the non-trivial flat line bundle satisfying 
\begin{equation}
c_1(L)=\alpha\in {\rm Tor}(H^2(M;\Z ))\cong {\rm Tor}(H_1(M;\Z ))
\label{torzio}
\end{equation}
and $a\in\ch^1_{L^2}(M)/2\pi\Lambda^1_{L^2}$ as above. This 
observation provides the description of $\widehat{M}$ in part (ii) of the 
lemma. $\Diamond$
\vspace{0.1in}	

\noindent Therefore in light of the previous lemmata and (\ref{szetszedes}) 
elements of $\ch^1_{L^2}(M)$ can be interpreted as flat $L^2$ perturbations 
of a connection $\nabla_L^0$ on a given line bundle $L$ over $M$. 

Now we can provide an explicit description of our relevant orbit space as 
follows.

\begin{lemma} The orbit space of gauge
inequivalent $L^2_l$ Abelian connections ($l\geqq 2$) over the trivial
line bundle $L_0$ with respect to the trivial flat connection
$\nabla^0_{L_0}$ and the Coulomb gauge condition admits a decomposition
\[\cb (\nabla^0_{L_0})\cong 
\frac{\ch^1_{L^2}(M)}{2\pi\Lambda^1_{L^2}}\times
\left( (\ch^1_{L^2}(M))^\perp\cap\overline{{\rm ker}\:\delta}\right)
,\:\:\:\:\:\overline{{\rm ker}\:\delta}\subset L^2_l(M;\Lambda^1M).\]
That is, it is a Hilbert space bundle over a finite dimensional torus
representing the gauge inequivalent flat connections on $L_0$ (cf. part  
{\rm (ii)} of Lemma \ref{laposlemma}).

We note again that if $M$ is compact then $L^2$ objects reduce to ordinary
de Rham cohomology.
\label{orbitterlemma}
\end{lemma}

\noindent {\it Proof.} Consider the trivial flat connection
$\nabla^0_{L_0}$ on $L_0$ and write any connection in the form
$\nabla_{L_0}=\nabla^0_{L_0}+a$ with $a\in\ca (\nabla^0_{L_0})$ as usual.
By the aid of Lemma \ref{coulomblemma} let us impose the Coulomb gauge
condition $\delta a=0$ on them.

An advantage of $L^2$ cohomology is that when $l\geqq 2$ the usual Hodge 
decomposition continues to hold \cite{hau-hun-maz} in the form
\begin{equation}
L^2_l(M;\Lambda^1M)\cong\overline{{\rm im}\:\dd}\oplus
\overline{{\rm im}\:\delta}\oplus\ch^1_{L^2}(M)
\label{hodge}
\end{equation}
and obviously $\overline{{\rm im}\:\delta}\oplus\ch^1_{L^2}(M)\subseteq
\overline{{\rm ker}\:\delta}$. If $\gamma\in\cu_{L_0}$ is any gauge
transformation on the trivial bundle then $\gamma^{-1}\dd\gamma$ is closed
moreover $\delta (a+\gamma^{-1}\dd\gamma )=0$ together with $\delta a=0$
implies that $\gamma^{-1}\dd\gamma$ is also co-closed and of course
$\ii\gamma^{-1}\dd\gamma\in L^2(M;\Lambda^1M)$. In other words
\[\ii\gamma^{-1}\dd\gamma\in\ch^1_{L^2}(M).\]
Hence Hodge decomposition shows that $\cu_{L_0}$ acts trivially on
$(\ch^1_{L^2}(M))^\perp\cap\overline{{\rm ker}\:\delta}$.
Consequently all elements here are gauge inequivalent while Lemma
\ref{laposlemma} shows that the finite dimensional torus
$\ch^1_{L^2}(M)/2\pi\Lambda^1_{L^2}\cong T^{b_1(M)}$ enumerates the
gauge inequivalent flat connections. $\Diamond$
\vspace{0.1in}

\noindent Now we move on to clarify the set of classical solutions used 
in this paper. This will also provide us with a physical interpretation of 
the space $\ch^2_{L^2}(M)$ as containing the curvatures of finite 
action classical solutions over a $4$-space $(M,g)$. First we introduce 
the class of four-manifolds considered in this paper. By definition this class 
contains all compact geometries moreover the so-called {\it asymptotically 
locally flat} (ALF) geometries. These latter sub-class contains 
non-compact but complete four dimensional Riemannian manifolds with a 
special asymptotical geometry. For the definition of an ALF space we refer 
to Sect. \ref{five}. Many important non-compact manifolds in 
mathematical physics are of ALF type. Examples are the multi-Taub--NUT 
spaces, the Riemannian Schwarzschild and Kerr geometries, etc.
It follows from \cite[Corollary 9]{hau-hun-maz} 
that an ALF space always has finite dimensional second $L^2$ cohomology. 

Writing $H^2(M;\R )\subseteq\ch^2_{L^2}(M)$ i.e., embedding 
compactly supported de Rham cohomology into $L^2$ cohomology by the 
unique harmonic representative in each cohomology class, we can suppose 
that $\Lambda^2\subset\ch^2_{L^2}(M)$ where $\Lambda^2$ is the 
integer lattice from (\ref{H2}). 

\begin{definition} Let $(M,g)$ be a compact or an ALF four-manifold. 
Consider its $2^{{\it nd}}$ $L^2$ cohomology. We define a co-compact lattice 
$\Lambda^2_{L^2}\subset\ch^2_{L^2}(M)$ as follows. We say that 
$\omega\in\Lambda^2_{L^2}$ if and only if

\begin{itemize}

\item[{\rm (i)}] either $\omega\in\Lambda^2$ i.e., it is a harmonic 
representative of a usual compactly supported integer de Rham cohomology class;

\item[{\rm (ii)}] or if $\omega\in\Lambda^2_{L^2}$ and $\omega\notin\Lambda^2$ 
then there exists an ${\rm U}(1)$-connection $\nabla^0_L$ on a 
line bundle $L$ with curvature $F_{\nabla^0_L}/2\pi\ii 
=\omega$ and this connection satisfies the strong holonomy condition at 
infinity.

\end{itemize}
\label{racsdefinicio} 
\end{definition}

\noindent For the definition of the strong holonomy condition we refer 
to \cite[Definition 2.1]{ete-jar}. Roughly speaking if a finite action 
${\rm U}(1)$-connection satisfies the strong holonomy condition then it 
has trivial holonomy at infinity hence in an appropriate $L^2_1$ norm it 
approaches the trivial flat ${\rm U}(1)$-connection on the infinitely distant 
boundary of $(M,g)$. This provides us that the action (\ref{hatas}), or 
equivalently (\ref{bovites}), of such a connection is an 
integer, cf. \cite[Theorem 2.2]{ete-jar}. Moreover it is clear that if 
$(M,g)$ is compact then part (ii) of Definition \ref{racsdefinicio} is vacuous.

Chern--Weil theory says that if $\omega\in\Lambda^2_{L^2}$ there 
exists a line bundle $L$ with $c_1(L)\in\Lambda^2\subseteq\Lambda^2_{L^2}$, 
the integer lattice in (\ref{H2}), and an ${\rm U}(1)$-connection 
$\nabla^0_L$ on the bundle such that $\omega =F_{\nabla^0_L}/2\pi\ii$. This 
connection is obviously a finite action classical solution to (\ref{maxwell}). 
Moreover it follows from part (ii) of Lemma \ref{laposlemma} and 
especially from (\ref{torzio}) that the remaining line bundles in 
(\ref{H2}) i.e., those with $c_1(L)\in{\rm Tor}(H^2(M;\Z ))$ also carry 
classical solutions $\nabla^0_L$ to (\ref{maxwell}) namely vacuum solutions 
i.e., flat connections. 

We conclude that all elements of the lattice 
$\Lambda^2_{L^2}\times{\rm Tor}(H^2(M;\Z ))$ give rise to 
classical solutions of the Abelian gauge theory (\ref{hatas}). 
In terms of (\ref{bovites}) or the projected fields 
$F^\pm_{\nabla^0_L}=\frac{1}{2}(F_{\nabla^0_L}\pm *F_{\nabla^0_L})$ all of 
them satisfy 
\begin{equation}
\frac{1}{8\pi^2}\int\limits_M\left( F_{\nabla^0_L}\wedge *F_{\nabla^0_L}
\pm F_{\nabla^0_L}\wedge F_{\nabla^0_L}\right)=
\frac{1}{4\pi^2}\int\limits_MF_{\nabla^0_L}^\pm\wedge *
F_{\nabla^0_L}^\pm\in\Z .
\label{kvantaltsag}
\end{equation}
We emphasize again that the validity of the energy quantization 
(\ref{kvantaltsag}) for the non-topological sector of $\Lambda^2_{L^2}$ 
follows from the boundary condition imposed on it in part (ii) of 
Definition \ref{racsdefinicio}.

\begin{remark}\rm In general on a given bundle 
these classical solutions are not unique because of two reasons. The first 
is if $b^2(M)<b^2_{L^2}(M)$ which can happen if $M$ is not compact, cf. 
Sect. \ref{four}. The second is if $\ch^1_{L^2}(M)\not=\{0\}$ which can 
happen if $M$ is not simply connected, cf. Lemma \ref{laposlemma}.
\end{remark}

\noindent Now we are in a position to carefully define the central object 
of our interest here. Consider a finite action classical solution 
$\nabla^0_L$ such that either its curvature represents an element in 
$\Lambda^2_{L^2}$ or it is a flat connection whose gauge class 
is in ${\rm Tor}(H^2(M;\Z ))$. In this case we write simply 
$[\nabla^0_L]\in\Lambda^2_{L^2}\times{\rm Tor}(H^2(M;\Z ))$ and these 
solutions will be referred to as {\it allowed classical solutions}. 
Moreover let $(M,g)$ be a compact or an ALF four-manifold. Our primary 
concern in this paper will be the calculation of the formal integral
\begin{equation}
Z(M,g,\tau ):=\sum\limits_{[\nabla^0_L]\in\left(\Lambda^2_{L^2}\times{\rm 
Tor}(H^2(M;\Z ))\right)}\frac{1}{{\rm Vol}(\cu_L)}
\int\limits_{\nabla_L\in\ca (\nabla^0_L)}{\rm e}^{-S(\nabla_L ,\tau )}
\DD\nabla_L
\label{particiofv}
\end{equation}
or equivalently
\[Z(M,g,\tau ):=\sum\limits_{[\nabla^0_L]\in\left(\Lambda^2_{L^2}
\times{\rm Tor}(H^2(M;\Z ))\right)}\:\:
\int\limits_{[\nabla_L]\in\cb (\nabla^0_L)} {\rm e}^{-S(\nabla_L ,\tau )}
\DD [\nabla_L]\]
which gives rise to the relevant part of the (Euclidean) {\it partition
function} of the supersymmetrized Abelian gauge theory (\ref{hatas})
over $(M,g)$. We note again that in the supersymmetric setting $\tau$ is
well-defined at the full quantum level hence its appearance in $Z(M,g ,\tau )$
is meaningful. Here $\DD\nabla_L$ denotes the
hypothetical (probably never definable) measure on the infinite
dimensional vector space $\ca (\nabla^0_L)$ and $\DD [\nabla_L]$ is the 
induced one on the orbit space $\cb (\nabla^0_L)$.

Notice that if $M$ is compact the summation reduces to a summation 
over line bundles i.e., over $H^2(M;\Z )$ in accord with \cite{wit1, 
wit2}. However if $M$ is non-compact then the summation is taken over more 
finite action classical solutions than bundles. Moreover we will see 
in Sect. \ref{four} that these allowed classical solutions are not the 
whole set of finite action classical solutions due to the strong 
holonomy condition from Definition \ref{racsdefinicio}.

Next we perform the straightforward summation in (\ref{particiofv}). 
Pick a bundle $L$ and an allowed classical solution $\nabla^0_L$ on it. 
With respect to $\nabla^0_L$ write a connection $\nabla_L\in\ca 
(\nabla^0_L)$ on the same bundle in the form $\nabla^0_L+a$. Consider the 
associated decomposition (\ref{szetszedes}) of the action. To be precise 
we regard $\nabla_L$ on $L\cong L\otimes L_0$ as a perturbation 
of $\nabla^0_L$ on $L$ with a connection $\dd +a$ on $L_0$. It then 
follows that the integral in (\ref{particiofv}) looks like
\[\int\limits_{[\nabla_L]\in\cb (\nabla^0_L)}{\rm e}^{-S(\nabla_L,\tau )}
\DD [\nabla_L]={\rm e}^{-S(\nabla^0_L,\tau )}
\int\limits_{[a]\in\cb (\nabla^0_L)}{\rm e}^{-\frac{{\rm Im}\tau}{8\pi}
\Vert\dd a\Vert^2_{L^2(M)}}\DD [a]\]
where $\DD [a]=\DD [\nabla_L]$ is the formal induced 
measure on $\cb (\nabla^0_L)$.

As we have seen bundles with $c_1(L)\in {\rm Tor}(H^2(M;\Z ))$ 
are flat i.e., the action (\ref{bovites}) vanishes along them, 
consequently the summation in (\ref{particiofv}) over the torsion part 
simply gives a numerical factor in the partition function 
leaving us with a summation over the free lattice part as follows:
\[Z(M,g,\tau )=\left\vert{\rm Tor}(H^2(M;\Z ))\right\vert
\left(\sum\limits_{[\nabla^0_L]\in\Lambda^2_{L^2}}\:
{\rm e}^{-S(\nabla^0_L,\tau )}\right)
\int\limits_{[a]\in\cb (\nabla^0_{L_0})}{\rm e}^{-\frac{{\rm 
Im}\tau}{8\pi}\Vert\dd a\Vert^2_{L^2(M)}}\DD [a]\]
where $\nabla^0_{L_0}$ is the trivial flat connection on the trivial line 
bundle $L_0$.

On the way we introduce our theta function. The left hand side of 
(\ref{kvantaltsag}) provides us with an 
$L^2$ quadratic form $q_M$ on $\Lambda^2_{L^2}\subset\ch^2_{L^2}(M)$ which 
is indefinite according to the splitting 
$\Lambda^2_{L^2}=\Lambda^+_{L^2}\times\Lambda^-_{L^2}$ into (anti)self-dual 
parts. The second $L^2$ Betti number also splits
like $b^2_{L^2}(M)=b^+_{L^2}(M)+b^-_{L^2}(M)$. Write $q_M=q_M^+\oplus q_M^-$ 
for the correspondig decomposition into definite parts. Note that if $M$ is 
compact and simply connected then $q_M$ is just the intersection form and 
$b^\pm_{L^2} (M)$ are just the usual signature decomposition of $b^2(M)$. Then 
$\vartheta_{q^\pm _M}:\C^+\rightarrow\C$ is defined by 
\[\vartheta_{q^\pm _M}(\tau ):=
\sum\limits_{n\:\in\:\underbrace{\Z\times\dots\times\Z}_{b^\pm_{L^2}(M)}}
{\rm e}^{\ii\pi q_M^\pm (n,n)\tau}\]
and has the following properties taking 
into account the unimodularity of the intersection form (cf. e.g. 
\cite[Sec. VII.6]{ser}): it is holomorphic on the upper half-plane 
moreover always satisfies the functional equations
\begin{equation}
\left\{\begin{array}{ll}
\vartheta_{q^\pm _M}(\tau +2)=\vartheta_{q^\pm _M}(\tau )& \mbox{(``level 
$2$ property'')}\\
\vartheta_{q^\pm _M}(-1/\tau )=(\tau/\ii )^{\frac{1}{2}b^\pm_{L^2}(M)}
\:\vartheta_{q^\pm _M}(\tau )& 
\mbox{(``modularity of weight $\frac{1}{2}b^\pm_{L^2}(M)$ property'')}
\end{array}\right.
\label{theta}
\end{equation}
where the square root is the principal value cut along the negative real 
axis.\footnote{In certain cases $\vartheta_{q^\pm _M}(\tau +1)=
\vartheta_{q^\pm _M}(\tau )$ also holds; for example if $M$ is a
compact spin manifold.}
By the aid of this function we proceed as follows. Making use of 
(\ref{kvantaltsag}) the summation in (\ref{particiofv}) over the remaining 
lattice $\Lambda^2_{L^2}$ gives
\begin{equation}
Z(M,g,\tau )=\left\vert {\rm Tor}(H^2(M;\Z ))\right\vert 
\vartheta_{q^+_M}(\tau )\vartheta_{q^-_M}(-\overline{\tau})
\int\limits_{\cb (\nabla^0_{L_0})}{\rm e}^{-\frac{{\rm 
Im}\tau}{8\pi}\Vert\dd a\Vert^2_{L^2(M)}}\:\DD [a].
\label{particiofv2}
\end{equation}
Note that if $b^-_{L^2}(M)>0$ then $Z(M,g, \tau )$ is not holomorphic in $\tau$.

We procced further and ask ourselves how to perform the remaining integral 
in (\ref{particiofv2}) by the aid of Lemma \ref{orbitterlemma}. It shows 
that integration in Coulomb gauge is to be taken over a finite dimensional 
torus and an infinite dimensional Hilbert space. The formal measure for the 
latter space is $\DD [a]=(\det '\triangle_0 )\DD a$ where 
$\det '\triangle_0$ is the formal determinant of the scalar Laplacian without 
its zero eigenvalues and $\DD a$ is some formal 
measure on the Hilbert space $L^2_l(M;\Lambda^1M)$ restricted to 
$\overline{{\rm ker}\:\delta}\subset L^2_l(M;\Lambda^1M)$. 
This determinant enters the story as the Faddeev--Popov determinant in 
Coulomb gauge in ${\rm U}(1)$ gauge theory. Of course $\det '\triangle_0$ is 
ill-defined; standard $\zeta$-function regularization might be used to define 
it. These issues will be investigated in the forthcoming sections.

Referring to Lemma \ref{orbitterlemma} since the action vanishes along 
the torus integrating over it we obtain 
\begin{equation}
\int\limits_{\cb (\nabla^0_{L_0})}{\rm e}^{-\frac{{\rm Im}\tau}{8\pi}
\Vert\dd a \Vert^2_{L^2(M)}}\:\DD [a]= 
{\rm Vol}\left(\frac{\ch^1_{L^2}(M)}{2\pi\Lambda^1_{L^2}}\right)
\det\: '\triangle_0
\left(\:\int\limits_{(\ch^1_{L^2}(M))^\perp\cap\overline{{\rm ker}\:\delta}}
{\rm e}^{-\frac{{\rm Im}\tau}{8\pi}\Vert\dd a\Vert^2_{L^2(M)}}\:\DD a\right) .
\label{particiofv2.5}
\end{equation}
In the Coulomb gauge, provided $a\in\ca (\nabla^0_{L_0})$ 
we can re-express the last term in (\ref{szetszedes}) as
\[\Vert\dd a\Vert^2_{L^2(M)}=(a\:,\:\triangle_1a)_{L^2(M)}\]
where $\triangle_1=\delta\dd +\dd\delta$ is the Laplacian acting on
$1$-forms and $(\cdot\: ,\:\cdot )_{L^2(M)}$ is the $L^2$ scalar product
on the space of $1$-forms.

Consequently, collecting all of our findings sofar, we obtain that we have 
eventually cut down the original integral (\ref{particiofv}) to a yet highly 
non-trivial formal integral
\begin{equation}
\int\limits_{(\ch^1_{L^2}(M))^\perp\cap\overline{{\rm ker}\:\delta}}
{\rm e}^{-\left(a\:,\:\frac{{\rm Im}\tau}{8\pi}\triangle_1a\right)_{L^2(M)}}
\:\DD a,\:\:\:\:\:\overline{{\rm ker}\:\delta}\subset L^2_l(M;\Lambda^1M).
\label{particiofv3}
\end{equation}
Our aim in the forthcoming sections will be to calculate this integral 
over various manifolds.

\begin{remark}\rm We make an important comment here related with 
the non-trivial moduli of flat connections. This comment is motivated by the 
difficulties we have to face when try to repeat this procedure in the 
non-compact case in Sect. \ref{four} and \ref{five}. As we already 
noted, given a 
fixed line bundle $L$ over $M$ with a fixed allowed classical solution 
$[\nabla^0_L]\in\Lambda^2_{L^2}\times{\rm Tor}(H^2(M;\Z ))$ to 
(\ref{maxwell}) then $\nabla^1_L=\nabla^0_L+a$ is another solution if 
$a\in\ch^1_{L^2}(M)$ since $a$ is simply a flat perturbation 
of $\nabla^0_L$. Therefore nothing prevents us to carry out the summation 
procedure again over $\Lambda^2_{L^2}\times{\rm Tor}(H^2(M,\Z ))$ but this 
time starting with the decomposition (\ref{szetszedes}) relative to 
$\nabla^1_L$ instead of $\nabla^0_L$. Fortunately since $S(\nabla^0_L,\tau )=
S(\nabla^1_L,\tau )$ we obtain the same result (\ref{particiofv2}) and 
(\ref{particiofv2.5}).

If $M$ is compact then all classical solutions arise this way hence the 
summation is unambigous consequently the shape of (\ref{particiofv2}) 
and (\ref{particiofv2.5}) is well-defined. However we will see that for 
non-compact manifolds the situation is not so simple.
\end{remark}


\section{Compact spaces}
\label{three}


In this section we calculate (\ref{particiofv3}) in the compact case via 
$\zeta$-function and heat kernel techniques. So throughout this section 
$(M,g)$ denotes a connected, compact, oriented Riemannian four-manifold 
without boundary. In this situation all the $L^2$ cohomology groups appeared 
sofar reduce to ordinary de Rham cohomology groups.

The operator $c\triangle_k$ with $c>0$ real constant is a 
positive symmetric operator on the orthogonal 
complement $(\ch^k(M))^\perp\subset L^2_l(M;\Lambda^kM)$ with $l\geqq 2$. 
By the finite dimensional analogue it is therefore natural to define 
the Gaussian-like integral to be
\[\int\limits_{(\ch^k(M))^\perp}{\rm e}^{-\left(a\:,\:
c\triangle_ka\right)_{L^2(M)}}\:\DD a
:=\pi^{\frac{1}{2}{\rm rk '}(c\triangle_k)}
\left(\det\:'\left(c\triangle_k\right)\right)^{-\frac{1}{2}}\]
where the regularized rank and the determinant is yet to be defined 
somehow. The familiar way to do this is by making use of 
$\zeta$-function regularization. Since the spectrum of the Laplacian over a 
compact manifold is non-negative real and discrete, one sets
\[\zeta_{\triangle_k}(s):=\sum\limits_{\lambda\in {\rm Spec}
\:\triangle_k-\{0\}}\lambda^{-s},\:\:\:\:\:\mbox{with $s\in\C$ and 
${\rm Re}\:s>0$ sufficiently large}\]
and observes that this function can be meromorphically continued over the 
whole complex plane (cf. e.g. \cite[Theroem 5.2]{ros}) having no pole at 
$s=0\in\C$. A formal calculation then convinces us that the regularized 
rank and the determinant of the Laplacian should be
\[{\rm rk'}\:\triangle_k:=\zeta_{\triangle_k}(0),\:\:\:\:\:\det\:'\triangle_k 
:={\rm e}^{-\zeta '_{\triangle_k}(0)}\]
yielding ${\rm rk'}(c\triangle_k)=\zeta_{\triangle_k}(0)$ and $\det '
(c\triangle_k)=c^{\zeta_{\triangle_k}(0)}
{\rm e}^{-\zeta '_{\triangle_k}(0)}$. Hence 
\begin{definition} Putting $c:={\rm Im}\tau /8\pi$ in the formal 
integral above we set
\[\int\limits_{(\ch^k(M))^\perp}
\!\!\!\!\!\!\!\!
{\rm e}^{-\left(a\:,\:
\frac{{\rm Im}\tau}{8\pi}\triangle_ka\right)_{L^2(M)}}\:\DD a:=
\pi^{\frac{1}{2}\zeta_{\triangle_k}(0)}\:{\rm e}^{\frac{1}{2}
\zeta '_{\triangle_k}(0)}
\left(\frac{{\rm Im}\:\tau}{8\pi}\right)^{-\frac{1}{2}\zeta_{\triangle_k}(0)}=
{\rm e}^{\frac{1}{2}\zeta '_{\triangle_k}(0)}
\left(\frac{{\rm Im}\:\tau}{8\pi^2}\right)^{-\frac{1}{2}
\zeta_{\triangle_k}(0)}.\]
Therefore the restricted formal integral (\ref{particiofv3}) taking 
place on the closed subspace $(\ch^1(M))^\perp\cap{\rm ker}\:\delta$ only is 
defined to be
\[\int\limits_{(\ch^1(M))^\perp\cap{\rm ker}\:\delta}
\!\!\!\!\!\!\!\!
{\rm e}^{-\left(a\:,\:\frac{{\rm 
Im}\tau}{8\pi}\triangle_1a\right)_{L^2(M)}}
\:\DD a:={\rm e}^{\frac{1}{2}\zeta '_{\triangle_1
\vert (\ch^1(M))^\perp\cap{\rm ker}\:\delta}(0)}
\left(\frac{{\rm Im}\:\tau}{8\pi^2}\right)^{-\frac{1}{2}\zeta_{\triangle_1
\vert (\ch^1(M))^\perp\cap{\rm ker}\:\delta}(0)}\]
where the restricted function $\zeta_{\triangle_1
\vert (\ch^1(M))^\perp\cap{\rm ker}\:\delta}$ is defined in the analogous way.
\label{integral.definiciok}
\end{definition}
\noindent Since this restricted $\zeta$-function 
is not easy to find we derive a sort of ``Fubini principle'' to 
obtain our integral successively from simpler ones.

\begin{lemma} With respect to Definition \ref{integral.definiciok} we find
\begin{equation}
\int\limits_{(\ch^1(M))^\perp\cap{\rm ker}\:\delta}{\rm e}^{-\left(
a\:,\:\frac{{\rm Im}\tau}{8\pi}\triangle_1a\right)_{L^2(M)}}\:\DD a =
{\rm e}^{\frac{1}{2}\zeta '_{\triangle_1}(0)-\zeta '_{\triangle_0}(0)}
\left(\frac{{\rm Im}\:\tau}{8\pi^2}\right)^{\frac{1}{2}
\left(\zeta_{\triangle_0}(0)-\zeta_{\triangle_1}(0)\right)}.
\label{keszvan}
\end{equation}
This is the value of the remaining integral (\ref{particiofv3}) in
the compact case.
\end{lemma}
 
\begin{remark}\rm Inserting (\ref{keszvan}) into (\ref{particiofv2.5}) 
and then into (\ref{particiofv2}) the calculation of the partition 
function (\ref{particiofv}) over compact spaces is now 
complete.\footnote{Including
the Faddeev--Popov determinant the full contribution
of the determinants to the partition function (\ref{particiofv}) is
$(\det '\triangle_1)^{-\frac{1}{2}}(\det '\triangle_0)^2$ that is,
${\rm e}^{\frac{1}{2}\zeta '_{\triangle_1}(0)-2\zeta '_{\triangle_0}(0)}$.
This could be further analyzed however we skip this here.}
\end{remark}

\noindent {\it Proof}. In addition to the Hodge decomposition (\ref{hodge}) we 
obviously know that ${\rm im}\:\delta\oplus\ch^1(M)\subseteq{\rm ker}\:\delta$ 
and ${\rm im}\:\dd\cap{\rm ker}\:\delta =\{ 0\}$ hence 
$L^2_l(M;\Lambda^1M)\cong {\rm im}\:\dd\oplus{\rm ker}\:\delta$. 
Intersecting this with $(\ch^1(M))^\perp$ and taking into account that 
${\rm im}\:\dd\cong (\ch^0(M))^\perp$ we obtain the further decomposition 
\begin{equation}
(\ch^1(M))^\perp\cong (\ch^0(M))^\perp\oplus\left( (\ch^1(M))^\perp 
\cap{\rm ker}\:\delta\right) .
\label{dekompozicio1}
\end{equation}
Applying (\ref{dekompozicio1}) in the compact case we can 
write any element $\ii a\in L^2_l(M;\Lambda^1M)$ 
uniquely in the form $a=\dd f +\alpha$ with $\ii f\in 
L^2_{l+1}(M;\Lambda^0M)$ a function and $\ii\alpha\in 
L^2_l(M;\Lambda^1M)$ satisfying $\delta\alpha =0$. A simple calculation 
ensures us that if $\ii f\in L^2_{l+2}(M;\Lambda^0M)$ then
\begin{equation}
\left( a\:,\:\triangle_1a\right)_{L^2(M)}=\left( 
\dd f +\alpha\:,\:\triangle_1(\dd f +\alpha)\right)_{L^2(M)}=
\left( f\:,\:\triangle^2_0f\right)_{L^2(M)}+
\left(\alpha\:,\:\triangle_1\alpha\right)_{L^2(M)}
\label{dekompozicio2}
\end{equation}
where $\triangle^2_0$ is the square of the scalar Laplacian on $(M,g)$.
Taking into account (\ref{dekompozicio1}) and 
(\ref{dekompozicio2}) we obtain that $({\rm Spec}\:\triangle_1-\{ 0\})= 
({\rm Spec}\:\triangle^2_0-\{0\})\sqcup ({\rm Spec}\:
\triangle_1\vert_{(\ch^1(M))^\perp\cap{\rm ker}\:\delta})$. This 
decomposition together 
with the proof of \cite[Theorem 5.2]{ros} ensures us that
\[\zeta_{\triangle_1}=\zeta_{\triangle^2_0}+
\zeta_{\triangle_1\vert(\ch^1(M))^\perp\cap{\rm ker}\:\delta}\]
consequently Definition \ref{integral.definiciok} yields the 
following Fubini-like successive formula
\[\int\limits_{(\ch^1(M))^\perp}
\!\!\!\!\!\!\!\!
{\rm e}^{-\left(a\:,\:\frac{{\rm Im}\tau}{8\pi}\triangle_1a\right)_{L^2(M)}}
\:\DD a=\left(\:\int\limits_{(\ch^0(M))^\perp}\!\!\!\!\!\!\!\!
{\rm e}^{-\left( f\:,\:\frac{{\rm Im}\tau}{8\pi}\triangle^2_0f\right)_{L^2(M)}}
\:\DD f\right)\!\!\!
\left(\:\int\limits_{(\ch^1(M))^\perp\cap{\rm ker}\:\delta}
\!\!\!\!\!\!\!\!\!\!\!\!\!\!
{\rm e}^{-\left(\alpha\:,\:\frac{{\rm Im}\tau}{8\pi}
\triangle_1\alpha\right)_{L^2(M)}}\:\DD\alpha\right) .\]
Taking into account 
that $\zeta_{\triangle^2_0}(s)=\zeta_{\triangle_0}(2s)$ 
hence ${\rm rk'}(c\triangle^2_0(0))=\zeta_{\triangle_0}(0)$ as well as 
$\det '(c\triangle^2_0)=c^{\zeta_{\triangle_0}(0)}
{\rm e}^{-2\zeta '_{\triangle_0}(0)}$ we obtain the lemma. $\Diamond$
\vspace{0.1in}

\noindent The $S$-duality properties of the partition function are 
concentrated 
in its $\vartheta_{q^+M}(\tau )\vartheta_{q^-M}(-\overline{\tau})$ 
term and the exponent of ${\rm Im}\tau/8\pi^2$ in (\ref{keszvan}). 
Hence we shall focus our attention to this exponent. Over a compact 
four-manifold $(M,g)$ without boundary it is well-known 
\cite[Theorem 5.2]{ros} that
\begin{equation}
\zeta_{\triangle_k}(0)=-\dim{\rm ker}\triangle_k+
\frac{1}{16\pi^2}\int\limits_M{\rm tr}(u^4_k)\dd V
\label{zetanull}
\end{equation}
where the sections $u^p_k\in C^\infty (M;{\rm End}(\Lambda^kM))$ with 
$p=0,1,\dots$ appear \cite [Chapter 3]{ros} in the coefficients of the 
short time asymptotic expansion of the heat kernel for the $k$-Laplacian
\[\sum\limits_{\lambda\in{\rm Spec}\triangle_k}
{\rm e}^{-\lambda t}\:\sim\:\frac{1}{(4\pi t)^2}\sum\limits_{p=0}^{+\infty}
\left(\:\:\int\limits_M{\rm tr}(u^p_k)\dd V\right) 
t^{\frac{p}{2}}\:\:\:\:\:\mbox{as $t\rightarrow 0$}.\]
These functions are expressible with the curvature of $(M,g)$ and one 
can demonstrate \cite[p. 340]{gil} that
\[u^4_0 =\frac{1}{360}\left( 2\vert R\vert^2-2\vert r\vert^2+5s^2\right)\]
and
\[{\rm tr}(u^4_1) =\frac{1}{360}\left( -22\vert R\vert^2+172\vert r\vert^2
-40s^2\right)\]
where $R$ is the Riemann, $r$ is the Ricci and $s$ is 
the scalar curvature of the metric on $M$. On substituting these into 
the exponent of ${\rm Im}\tau /8\pi^2$ in (\ref{keszvan}) we find that it 
takes the shape
\[\frac{1}{2}\left(\zeta_{\triangle_0}(0)-\zeta_{\triangle_1}(0)\right)=
\frac{1}{2}\left( b^1(M)-b^0(M)+\frac{1}{\pi^2}\int\limits_M
\left(\frac{1}{120}\vert R\vert^2-\frac{87}{2880}\vert r\vert^2+
\frac{1}{128}s^2\right)\dd V\right) .\]

The time has come to write down the modular weights of the partition 
function (\ref{particiofv}). First note that
\[{\rm Im}\left(-\frac{1}{\tau}\right)=
\frac{1}{\tau\overline{\tau}}{\rm Im}\:\tau\]
that is, it is modular of holomorphic and anti-holomorphic weights 
$(-1,-1)$ respectively. Secondly it follows from (\ref{theta}) that up 
to $\ii$'s the modular weight of $\vartheta_{q_M^+}(\tau )$ and 
$\vartheta_{q_M^-}(-\overline{\tau})$  is $\frac{1}{2}b^\pm_{L^2}(M)$ hence 
referring to the shape of the partition function in (\ref{particiofv2}) we 
find that its holomorphic and anti-holomorphic modular weights 
$(\alpha ,\beta )$ from (\ref{eredmeny}) are
\begin{equation}
\frac{1}{4}\left(\chi (M)\pm\sigma (M)-\frac{1}{\pi^2}\int\limits_M
\left( \frac{1}{60}\vert R\vert^2-\frac{87}{1440}\vert r\vert^2
+\frac{1}{64}s^2\right)\dd V\right)
\label{sulyok}
\end{equation}
(to be precise the ``$+$'' is for holomorphic and the ``$-$'' is for 
anti-holomorphic). In this formula $\chi(M)\pm\sigma 
(M)=2b^0(M)-2b^1(M)+2b^\pm (M)$ are the linear combinations of the Euler 
characteristic and the signature of the manifold. 

Since (\ref{sulyok}) is not zero in general we 
conclude that in the naive theory (\ref{hatas}) $S$-duality 
breaks down because its partition function (\ref{particiofv}) is 
not modular in $\tau$. Moreover, comparing (\ref{sulyok}) with Witten's 
calculation \cite{wit1, wit2} whose result is simply the topological 
term $\frac{1}{4}(\chi (M)\pm\sigma (M))$ we find an analytic correction 
which vanishes only in the rare situation if $(M,g)$ 
happens to be flat (there exist only $27$ connected compact orientable 
flat four-manifolds and as many as $74$ if the non-orientable ones are 
also included \cite{hil}).

Nevertheless this curvature correction does not destroy the main 
conclusion in \cite{wit1,wit2} namely that Abelian $S$-duality over a 
compact space can be restored within the framework of local 
quantum field theories by adding appropriate gravitational terms to the 
Lagrangian (\ref{hatas}) in order to cancel the modular anomaly coming from 
(\ref{sulyok}). As it was observed by Witten \cite{wit1, wit2} it is 
quite remarkable that although the individual Betti numbers not, their 
combinations $\chi (M)\pm\sigma (M)$ are expressible as 
integrals of local curvature densities (cf. e.g. \cite[p. 
370-371]{bes}): the Gauss--Bonnet--Chern theorem gives
\begin{equation}
\chi (M)= \frac{1}{8\pi^2}\int\limits_M{\rm tr} (R\wedge 
*R)=\frac{1}{8\pi^2}\int\limits_M\left(\vert R\vert^2-
\left\vert r-\frac{s}{4}g\right\vert^2\right)\dd V
\label{euler}
\end{equation}
for the Euler characteristic and the Hirzebruch signature theorem asserts that
\begin{equation}
\sigma (M)=-\frac{1}{24\pi^2}\int\limits_M{\rm tr}(R\wedge R)
\label{szignatura}
\end{equation}
holds for the signature. 
Hence because obviously the full weights (\ref{sulyok}) 
continue to be integrals of local densities, modular anomaly cancels by 
adding further gravitational terms to the Lagrangian 
(\ref{hatas}) (called ``$c$-numbers'') however they are not of the form 
$a(\tau ,\overline{\tau} ){\rm tr}(R\wedge *R)+b(\tau ,
\overline{\tau}){\rm tr}(R\wedge R)$ as claimed in \cite{wit1,wit2}. 
Their shape can be read off from (\ref{sulyok}). 


\section{Significance of the strong holonomy condition}
\label{four}


In this section we clarify the role of the strong holonomy condition 
imposed on connections in Definition \ref{racsdefinicio}. This 
condition excludes certain finite action classical solutions from 
the set of connections contributing to the partition function 
(\ref{particiofv}). It will turn out now that without this condition the 
partition function would have pathological behaviour. To show this we take 
the underlying space to be the multi-Taub--NUT spaces; these are quite 
important hyper-K\"ahler ALF spaces (see next section).

So let $(M_V,g_V)$ be the $1$-Taub--NUT space. This is a non-flat 
hyper-K\"ahler geometry on $\R^4\cong M_V$ and represents the $s=1$ 
member of the multi-Taub--NUT series ($s\in\N$ refers to 
the number of NUTs). For our purposes here we refer to \cite{ete-jar} 
for a description of this space. First of all one can demonstrate that 
it is an ALF space and its unique (up to scale) non-topological $L^2$ 
harmonic $2$-form arises as follows. Put an orientation onto $M_V$ induced by 
any complex structure in the hyper-K\"ahler family. As it is well-known the 
$1$-Taub--NUT space admits a non-trivial $L^2$ harmonic $2$-form $\omega$. 
This $2$-form can be constructed as the exterior derivative of the metric 
dual of the Killing field generating an isometric action of ${\rm U}(1)$ on the 
$1$-Taub--NUT space \cite{gib}. One can also obtain it by the conformal 
rescaling method when hunting for ${\rm SU}(2)$ anti-instantons 
\cite{ete-hau2}. Taking into account that $H^2(M_V ;\R )=\{0\}$ the condition 
$\dd\omega =0$ is equivalent to the existence of an imaginary valued $1$-form 
$A$ such that $\dd A=\ii\omega$. In other words there exists an 
${\rm U}(1)$-connection $\nabla^1_{L_0}:=\dd +A$ on the trivial 
bundle $L_0\cong M_V\times\C$ such that $\ii F_{\nabla^1_{L_0}}\in 
L^2(M_V;\Lambda^2M_V)$ as well as $\dd F_{\nabla^1_{L_0}}=0$ and 
$\delta F_{\nabla^1_{L_0}}=0$. That is, $\nabla^1_{L_0}$ is a finite 
action classical solution to the Maxwell equations (\ref{maxwell}) over 
the $1$-Taub--NUT space (moreover it is anti-self-dual). As a consequence of 
the linearity of the Abelian gauge theory and that $\omega$ is topologically 
trivial, for all $c\in\R$ the rescaled connection $\nabla^c_{L_0}:=\dd +cA$ 
is another gauge inequivalent solution on the same bundle with action 
proportional to $c^2$. In fact this $1$-parameter family of solutions is 
anti-self-dual hence from (\ref{bovites}) we obtain that
\[S(\nabla^c_{L_0},\tau )=\ii\pi\overline{\tau}c^2,\:\:\:\:\:c\in\R\]
in contrast to the quantized nature of (\ref{kvantaltsag}). Now picking 
any $L^2_l$ perturbation $a\in\ca (\nabla^0_{L_0})$ as before we find via 
(\ref{szetszedes}) that
\[S(\nabla^c_{L_0}+a,\tau )=\ii\pi\overline{\tau}c^2 +\frac{{\rm 
Im}\:\tau}{8\pi}\Vert\dd a\Vert^2_{L^2(M_V)}.\]
Moreover note that if $c\not= 0$ then of course 
$\nabla^c_{L_0}\notin\ca (\nabla^0_{L_0})$ (otherwise 
we would find $F_{\nabla^c_{L_0}}=0$ by Stokes' theorem) hence 
treating these solutions not as $L^2_l$ perturbations of the trivial 
flat connection is correct even from the functional analytic viewpoint. 

Suppose now that we want to calculate the partition function by simply 
mimicing the calculation in \cite{wit1, wit2} designed for the compact 
case i.e., integrating over connections on a given bundle and then summing 
over line bundles. Then, taking into account that $M_V\cong\R^4$ is simply 
connected hence there are no non-trivial flat connections 
and there is only one line bundle on it the partition function 
(\ref{particiofv}) is expected to look like
\begin{equation}
Z(M_V,g_V,\tau )=\left(\:\int\limits_{-\infty}^{+\infty}
{\rm e}^{\ii\pi (-\overline{\tau})c^2}\dd c\right)\det\:'\triangle_0
\left(\:\:\int\limits_{(\ch^1_{L^2}(M))^\perp\cap
\overline{{\rm ker}\:\delta}}{\rm e}^{-\left( a\:,\:
\frac{{\rm Im}\tau}{8\pi}\triangle_1a\right)_{L^2(M)}}\:\DD a\right) .
\label{rossz}
\end{equation}
This formula replaces (\ref{particiofv2}) and the last term is 
just the formal integral already appeared in (\ref{particiofv3}). Since 
$\dd c$ is {\it a fortiori} the Lebesgue measure on $\R$ and 
${\rm Im}\:\tau >0$ the Gaussian integral in the front converges and we 
plainly obtain
\[Z(M_V,g_V,\tau )=
(\ii /\overline{\tau})^{\frac{1}{2}}\:(\det\:'\triangle_0)
\left(\:\:\int\limits_{(\ch^1_{L^2}(M_V))^\perp\cap\:
\overline{{\rm ker}\:\delta}}{\rm e}^{-\left( a\:,\:
\frac{{\rm Im}\tau}{8\pi}\triangle_1a\right)_{L^2(M_V)}}\:\DD a\right) .\]
However this formula is far from being able to satisfy something like 
(\ref{eredmeny}) moreover the expected weights of 
$\tau$ and $\overline{\tau}$ do not look like in the compact case 
consequently modular anomaly cannot be cancelled with the same mechanism.

The situation gets even worse over the multi-Taub--NUT 
spaces $(M_V, g_V)$ with $s>1$ NUTs. In this case the aforementioned 
anti-self-dual $L^2$ solutions $\nabla^c_{L_0} =\dd +cA$ for all $c\in\R$ 
exist on $L_0$ (cf. \cite{ete-jar}). Hence if $L$ is a generic line 
bundle over $M_V$ carrying a classical 
solution $\nabla^0_L$ of the Maxwell equations then we obtain a similar 
$1$-parameter family $\nabla^c_L:=\nabla^0_L+cA$ on $L\cong L\otimes L_0$. 
Consequently in this case we would end up with more complicated 
ill expressions for the partition function.

Observe that we have run into this divergence problem not 
because of using some inappropriate regularization method for the infinite 
dimensional integral (\ref{particiofv3}); rather the problem arose from 
the multi-Taub--NUT geometry itself in the sense that it possesses 
``too many'' finite action classical solutions, enough to distort 
the partition function. We encounter the same difficulty over the 
Riemannian Schwarzschild and Kerr geometries, too (these are also 
ALF spaces, cf. \cite{ete-hau1}). These problematic 
solutions are truely non-compact phenomena in the sense that the $L^2$ 
harmonic $2$-forms their curvatures represent are non-topological. 
This is the moment where the powerful nature of Definition 
\ref{racsdefinicio} shows up: it rules out most of these pathological 
solutions but not all of them! We will see in Sect. \ref{five} that if all of 
these non-topological problematic solutions are excluded then the modular 
properties would be destroyed again but in a different way. 

In fact the remaining classical finite action 
solutions excluded by Definition \ref{racsdefinicio} 
represent surface operators (cf. e.g. \cite{tan1, tan2}) attached to the 
infinitely distant surface $i:B_{\infty}\subset X$ (see 
Sect. \ref{five}) where $X$ is a natural compactification of the 
ALF space \cite{ete, ete-jar, hau-hun-maz}. That is, they provide us with 
observables of the form $\nabla^c_{\tilde{L}}\mapsto 
O_{B_{\infty}}(\nabla^c_{\tilde{L}}):=
\exp (\int_{B_{\infty}}i^*F_{\nabla^c_{\tilde{L}}})\in\C$ in our theory where 
$\nabla^c_{\tilde{L}}$ is a singular connection on the 
extended line bundle $\tilde{L}$ over $X$. 

We also note that if the Gaussian integral in (\ref{rossz}) could be somehow 
replaced by
\[\frac{1}{\ii}\int\limits_{-\infty}^{+\infty}
{\rm e}^{\ii\pi(-\overline{\tau})c^2}\cot (\pi c)\dd c
:=\lim\limits_{\varepsilon\rightarrow 0}\frac{1}{\ii}
\int\limits_{-\infty}^{+\infty}{\rm e}^{\ii\pi 
(-\overline{\tau})(c+\ii\varepsilon )^2}\cot (\pi (c+\ii\varepsilon ))\dd c
=\vartheta_{q^-_{M_V}}(-\overline{\tau})\]
i.e., the restricted Feynman measure would turn out to be the singular measure 
$\frac{\cot (\pi c)}{\ii}\dd c$ then modular properties of $Z(M_V,g_V,\tau )$ 
are also recovered. 

Accepting the definition of the partition function via 
(\ref{particiofv}) a second problem arises if one tries to calculate 
the infinite dimensional integral (\ref{particiofv3}) over an 
ALF space akin to (\ref{keszvan}) i.e., by making use of a 
$\zeta$-function regularization. This time one has to face the 
problem that in general the spectra of Laplacians on $k$-forms are not 
discrete hence the existence of their $\zeta$-functions is not obvious. 
This obstacle will be resolved in the next section by a straightforward 
regularization method based on truncating the non-compact space and 
imposing Dirichlet boundary condition on the boundary.


\section{Asymptotically locally flat spaces}
\label{five}


In this section we repeat our calculations of (\ref{particiofv}) over 
ALF spaces. The calculation goes along the same lines as in the compact space 
with obvious technical modifications.

ALF spaces are also referred to sometimes as {\it gravitational instantons of 
ALF type} in the broad or narrow sense if in addition their metric is 
Ricci-flat or hyper-K\"ahler respectively. Ricci flat examples are the 
Riemannian Schwarzschild and Kerr manifolds while the flat space 
$\R^3\times S^1$ and the multi-Taub--NUT spaces (also called $A_k$ ALF or ALF 
Gibbons--Hawking spaces) and the $D_k$ ALF spaces (including the 
Atiyah--Hitchin manifold as the $D_0$ case) provide hyper-K\"ahler 
examples. The multi-Taub--NUT spaces apparently are subject to recent 
investigations in mathematical physics \cite{che1, che2, ete-sza, wit3}.

To begin with, we recall the definition of these spaces taken 
from \cite{ete-jar} for instance. Let $(M,g)$ be a connected, oriented 
Riemannian four-manifold. This
space is called an {\it asymptotically locally flat (ALF) space} if
the following holds. There is a compact subset $K\subset M$ such
that $M\setminus K =W$ and $W\cong N\times\R^+$, with $N$ being a
connected, compact, oriented three-manifold without boundary
admitting a smooth $S^1$-fibration
\[\pi :N\stackrel{F}{\longrightarrow}B_\infty\]
whose base space is a compact Riemann surface $B_\infty$.
For the smooth, complete Riemannian metric $g$ there exists a
diffeomorphism $\phi :N\times\R^+\rightarrow W$ such that
\[\phi^*(g\vert_W)=\dd\rho^2+\rho^2(\pi^*g_{B_\infty})'+h'_F\]
where $g_{B_\infty}$ is a smooth metric on $B_\infty$, $h_F$ is a
symmetric 2-tensor on $N$ which restricts to a metric along the
fibers $F\cong S^1$ and $(\pi^*g_{B_\infty})'$ as well as $h'_F$ are
some finite, bounded, smooth extensions of $\pi^*g_{B_\infty}$ and $h_F$
over $W$, respectively. That is, we require $(\pi^*g_{B_\infty})'(\rho )\sim 
O(1)$ and $h'_F(\rho )\sim O(1)$ and the extensions for $\rho <+\infty$ 
preserve the properties of the original fields. Furthermore,
we require the Riemann curvature $R$ of $g$ to decay like
\begin{equation}
\vert\phi^*(\nabla^k R\vert_W)\vert\sim
O(\rho^{-3-k}),\:\:\:\:\:k=0,1,2,\dots
\label{gorbulet}
\end{equation}
where $R$ is regarded as a map $R:C^\infty (M;\Lambda^2M)\rightarrow 
C^\infty (M;\Lambda^2M)$ and its pointwise norm is calculated accordingly 
in an orthonormal frame. The definition of the metric shows that the 
volume of our spaces is  infinite however from the curvature decay 
(\ref{gorbulet}) follows that both $\chi_{L^2}(M)$ and $\sigma_{L^2}(M)$ 
defined by (\ref{euler}) and (\ref{szignatura}) respectively, remain finite.

Take any ALF space $(M,g)$ as above. Fix a real number $\rho >0$ and let 
$(\overline{M}_\rho ,g\vert_{\overline{M}_\rho})$ be the truncated 
manifold i.e., $\overline{M}_\rho$ is a four-manifold with 
connected boundary $\partial\overline{M}_\rho$ such that the distance of any 
point $x\in\overline{M}_\rho$ satisfies $d(x_0,x)\leqq\rho$ with 
respect to an interior point $x_0\in M_\rho$. 

With $l\geqq 2$ define the spaces 
\[LD^2_l(\overline{M}_\rho ;\Lambda^k\overline{M}_\rho):=
\{\mbox{$\varphi\in L^2_l(\overline{M}_\rho ;\Lambda^k\overline{M}_\rho )
\:\vert\:\varphi\vert_{\partial\overline{M}_\rho}=0$ and 
$\delta\varphi\vert_{\partial\overline{M}_\rho}=0$ a.e.}\}\]
and 
\[\ch_D^k(\overline{M}_\rho ):=
\{\mbox{$\varphi\in LD^2_l(\overline{M}_\rho ;\Lambda^k\overline{M}_\rho)
\:\vert\:\dd\varphi =0$ and $\delta\varphi =0$}\}
\cong{\rm ker}(\triangle_k\vert_{\overline{M}_\rho}).\]
Also write $b_D^k(\overline{M}_\rho):=\dim\ch_D^k(\overline{M}_\rho )$
for the corresponding Dirichlet--Betti numbers. Note that pushing 
$\partial\overline{M}_\rho$ toward infinity these numbers satisfy 
\begin{equation}
b_D^k(\overline{M}_\rho)\leqq b^k_{L^2}(M)
\label{betti}
\end{equation}
where these latter numbers are 
the true $L^2$ Betti numbers of the original open manifold. 

Now we are in a position to calculate the partition function over an 
ALF space as follows. First define $Z(M,g,\tau )$ by 
(\ref{particiofv}) including the remaining integral (\ref{particiofv3}). 
To calculate this integral we observe the following things. Let us work over 
$(\overline{M}_\rho ,g\vert_{\overline{M}_\rho})$. Referring 
to the Hodge decomposition theorem with respect to the Dirichlet boundary 
condition 
\cite[Proposition 5.9.8]{tay} the decomposition (\ref{dekompozicio1}) 
continues to hold providing us with the validity of (\ref{dekompozicio2}). 
Secondly, under the Dirichlet condition the 
Laplacian has a discrete spectrum on $LD^2_l$ consequently introducing   
$\zeta_{\triangle_k}$ also makes sense. This enables us to express the 
integral (\ref{particiofv3}) over 
$(\overline{M}_\rho,g\vert_{\overline{M}_\rho})$ with (\ref{keszvan}). 
Then sending the boundary to infinity if the limit exists then the 
integral (\ref{particiofv3}) is 
defined by this limit. In this way we also obtain an expression for the 
original partition function (\ref{particiofv}).
 
Now following \cite{gil} we quickly summarize the changes 
of the heat kernel formul{\ae} in the situation of the Dirichlet boundary 
condition. The expression for the zero value of the $\zeta$-function looks like
\[\zeta_{\triangle_k\vert_{\overline{M}_\rho}}(0)=-\dim{\rm ker}
(\triangle_k\vert_{\overline{M}_\rho})+\frac{1}{16\pi^2}
\int\limits_{\overline{M}_\rho}{\rm tr}(u^4_k)\dd V+\frac{1}{16\pi^2}
\int\limits_{\partial\overline{M}_\rho}{\rm tr}(v^4_k)
\dd V\vert_{\partial\overline{M}_\rho}\]
and the integrals are again defined via the short time heat kernel 
expansion
\[\sum\limits_{\lambda\in{\rm Spec}(\triangle_k\vert_{\overline{M}_\rho})}
{\rm e}^{-\lambda t}\:\sim\:\frac{1}{(4\pi t)^2}\sum\limits_{p=0}^{+\infty}
\left(\:\:\int\limits_{\overline{M}_\rho}{\rm tr}(u^p_k)\dd V+
\int\limits_{\partial\overline{M}_\rho}{\rm tr}(v^p_k)
\dd V\vert_{\partial\overline{M}_\rho}\right)  
t^{\frac{p}{2}}\:\:\:\:\:\mbox{as $t\rightarrow 0$}. \]
The curvature expressions on the bulk are similar to the compact case 
\cite[Theorem 4.5.1]{gil}:
\[u^4_0=\frac{1}{360}\left( 2\vert R\vert^2
-2\vert r\vert^2+5s^2+12\triangle_0s\right)\]
and 
\[{\rm tr}(u^4_1)=\frac{1}{360}\left( 
-22\vert R\vert^2+172\vert r\vert^2-40s^2+48\triangle_0s\right) .\]
However the contribution of the boundary is quite complicated. Consider a 
collar $U\subset\overline{M}_\rho$ of the boundary which looks like 
$U\cong \partial\overline{M}_\rho\times (-1,0]$. Pick an 
orthonormal frame field $(e_1,\dots,e_4)$ along $U$ such that $e_4$ is 
orthogonal to the boundary. Then define the second fundamental form of the 
boundary by 
\[\Pi_{ij}:=g\vert_{\overline{M}_\rho}(\nabla_{e_i}e_j\:,\:e_4),
\:\:\:\:\:i,j=1,2,3\]
as well as let $\overline{\nabla}$ be the restriction of the 
Levi--Civita connection to the boundary. The boundary contributions are 
\cite[Theorem 4.5.1]{gil}
\begin{eqnarray}
v^4_0 & =& \frac{1}{360}\Big( -138\nabla_{e_4}s +140s\Pi_{ii} 
+4R_{i4i4}\Pi_{jj}-12R_{i4j4}\Pi_{ij}+4R_{ijkj}\Pi_{ik}\nonumber\\
& + & 24\overline{\nabla}^2_{e_i,e_i}\Pi_{jj}
+\frac{40}{21}\Pi_{ii}\Pi_{jj}\Pi_{kk}-\frac{88}{7}\Pi_{ij}\Pi_{ij}\Pi_{kk}
+\frac{320}{21}\Pi_{ij}\Pi_{jk}\Pi_{ik}\Big)\nonumber 
\end{eqnarray}
and
\begin{eqnarray}
{\rm tr}(v^4_1) & =& \frac{1}{360}\Big( -192\nabla_{e_4}s+200s\Pi_{ii}
+16R_{i4i4}\Pi_{jj}-48R_{i4j4}\Pi_{ij}+16R_{ijkj}\Pi_{ik}\nonumber\\
& +& 96\overline{\nabla}^2_{e_i,e_i}\Pi_{jj}
+\frac{160}{21}\Pi_{ii}\Pi_{jj}\Pi_{kk}-\frac{352}{7}\Pi_{ij}\Pi_{ij}\Pi_{kk}
+\frac{1280}{21}\Pi_{ij}\Pi_{jk}\Pi_{ik}\Big) .\nonumber
\end{eqnarray}
Now we demonstrate that all contributions
from the boundary get vanish as we move the boundary toward infinity.
\begin{lemma}
Let $(M,g)$ be an ALF space carefully defined above and let 
$(\overline{M}_\rho, g\vert_{\overline{M}_\rho} )$ be its truncation. 

Then
\[\lim\limits_{\rho\rightarrow +\infty}
\int\limits_{\overline{M}_\rho}(\triangle_0s)\dd V=0\:\:\:\:\:{\it and}
\:\:\:\:\:\lim\limits_{\rho\rightarrow +\infty}
\int\limits_{\partial\overline{M}_\rho}(\nabla_{e_4}s)
\dd V\vert_{\partial\overline{M}_\rho}=0\]
moreover $\Pi_{ij}\sim O(\rho^{-1})$ for all $i,j=1,2,3$.
\end{lemma}

\noindent{\it Proof.} Referring to the definition of an ALF space we find 
that ${\rm Vol}(\partial\overline{M}_\rho )\sim O(\rho^2)$ and from the 
curvature decay (\ref{gorbulet}) we get 
$\vert\nabla r\vert\sim O(\rho^{-4})$ and 
$\vert\nabla s\vert\sim O(\rho^{-4})$. Regarding the first integral, applying 
Stokes' theorem to convert it into an integral of $\nabla r$ along the 
boundary, the result follows. The second integral also decays in a 
straightforward way. 

Finally, since asymptotically $\Pi_{ij}(x)\sim\Gamma^4_{ij}(x)$ for all 
$x\in\partial\overline{M}_\rho$, a simple calculation yields 
that over an ALF space $\Pi$ also decays as claimed. $\Diamond$
\vspace{0.1in}

\noindent We proceed further and take the limit $\rho\rightarrow 
+\infty$ to recover the original ALF space $(M,g)$. Applying the lemma 
and the curvature decay (\ref{gorbulet}) we can see that all the terms 
involving $\triangle_0s$, $\nabla_{e_4}s$ and $\Pi$ do not 
contribute to the integrals. Therefore, since $\dim{\rm 
ker}(\triangle_k\vert_{\overline{M}_\rho})=b^k_D(\overline{M}_\rho)$ 
we obtain over the original ALF space $(M,g)$ that 
the exponent of ${\rm Im}\tau /8\pi^2$ in (\ref{keszvan}) converges 
and gives again
\[\frac{1}{2}\left( b^1_D(M)-b^0_D(M)+\frac{1}{\pi^2}\int\limits_M
\left(\frac{1}{120}\vert R\vert^2-\frac{87}{2880}\vert r\vert^2+   
\frac{1}{128}s^2\right)\dd V\right) <+\infty\]
as expected. Note that by 
(\ref{betti}) the limiting Dirichlet--Betti numbers $b^k_D(M)$ are well-defined.

From (\ref{particiofv2}) we already know that the contributions of the 
$\vartheta$-functions to the modular weights are $b^{\pm}_{L^2}(M)$ 
hence we eventually obtain that over an ALF space $(M,g)$ 
if the partition function (\ref{particiofv}) exists\footnote{We continue to 
abandon questions about the existence of the contribution of the 
determinants to the partition function in the limiting ALF case i.e., we do 
not check the limit of 
${\rm exp}\left(\frac{1}{2}\zeta '_{\triangle_1\vert_{\overline{M}_\rho}}(0)-
2\zeta '_{\triangle_0\vert_{\overline{M}_\rho}}(0)\right)$ as
$\rho\rightarrow +\infty$.} then its modular weights are equal to:
\begin{equation}
\frac{1}{4}\left( 2b^0_D(M)-2b^1_D(M)+2b^\pm _{L^2}(M)
-\frac{1}{\pi^2}\int\limits_M\left(\frac{1}{60}\vert R\vert^2
-\frac{87}{1440}\vert r\vert^2+\frac{1}{64}s^2\right)\dd V\right) .
\label{alfsulyok}
\end{equation}
It is worth comparing this with the compact case (\ref{sulyok}). The 
maximum principle provides us that $b^0_D(\overline{M}_\rho )=0$ 
hence $b^0_D(M)=0$ too, moreover Yau's theorem \cite{yau} gives 
$b^0_{L^2}(M)=0$ since an ALF space is complete and has infinite volume. 
Additionally, since \cite[Proposition 5.9.9]{tay} 
\[b^k_D(\overline{M}_\rho )=\dim H^k(\overline{M}_\rho 
,\partial\overline{M}_\rho\: ;\R )\] 
the part of the relative de Rham exact sequence
\[\begin{array}{ll} 
& \{ 0\}\cong H^0(\overline{M}_\rho ,\partial\overline{M}_\rho\: ;\R )
\longrightarrow H^0(\overline{M}_\rho\: ;\R )\longrightarrow 
H^0(\partial\overline{M}_\rho\: ;\R )\\
         &   \\
& \longrightarrow H^1(\overline{M}_\rho,\partial\overline{M}_\rho\: ;\R )
\longrightarrow H^1(\overline{M}_\rho\: ;\R )
\longrightarrow H^1(\partial\overline{M}_\rho\: ;\R )
\end{array}\]
shows that if $H^1(\partial\overline{M}_\rho\: ;\R )=\{ 0\}$ then 
$b^1_D(M)=b^1(M)$ moreover $b^1(M)=b^1_{L^2}(M)$ by the remark of Lemma 
\ref{laposlemma}. Consequently for an ALF space 
$(M,g)$ with $M=K\cup W$ and $W$ representing its neck, if the mild 
topological condition 
\begin{equation}
H^1 (W;\R )=\{ 0\}
\label{perem}
\end{equation}
holds for its neck then 
\[2b^0_D(M)-2b^1_D(M)+2b^\pm_{L^2}(M)=\chi_{L^2}(M)\pm\sigma_{L^2}(M)\]
that is, the modular weights can be written as integrals of the 
curvature like in the case of compact spaces. Consequently the modular anomaly 
cancels as before (cf. Sect. \ref{three}). For example, the 
multi-Taub--NUT family (whose members are 
hyper-K\"ahler ALF spaces) satisfies (\ref{perem}) while the Riemannian 
Schwarzschild and Kerr spaces (which are Ricci flat ALF spaces) not.


\section{Conclusion}


Let us summarize our approach to the partition function of 
Abelian gauge theory over ALF geometries. 

We define the partition function via (\ref{particiofv}) which is a natural 
generalization of Witten's formula for compact manifolds. The crucial concept 
here is to impose the strong holonomy condition at infinity to attenuate the 
finite action classical solutions contributing to the partition function.
After summation the resulting integral (\ref{particiofv3}) is 
calculated by a regularization method based on truncating the original 
manifold and imposing the Dirichlet boundary condition on the boundary 
and then sending the boundary to infinity.

As a result the partition function can be calculated and up to a 
numerical factor it is a level $2$ modular form like in the compact 
case i.e., it transforms as (\ref{eredmeny}) with $b^\pm (M)$ replaced 
by $b^\pm _{L^2}(M)$. The (anti-)holomorphic modular weights are in 
general not zero hence $S$-duality fails in the naive theory (\ref{hatas}). 
However, if in addition 
the simple condition (\ref{perem}) on the topology of the infinity of the 
ALF space holds then the modular anomaly can be removed by 
adding gravitational counter terms to the Lagrangian. Consequently $S$-duality 
can be restored within the realm of local quantum field theories as in 
the compact case. 

This program can be carried out at least in the case 
of the multi-Taub--NUT geometries. For instance for the $1$-Taub--NUT 
space $(M_V,g_V)$ with its standard orientation the result is as follows. 
On $M_V\cong\R^4$ there are no non-trivial flat connections and regarding  
$L^2$ cohomology one knows that $b^0_{L^2}(M_V)=b^1_{L^2}(M_V)=0$ 
moreover $b^2_{L^2}(M_V)=b^-_{L^2}(M_V)=1$ 
yielding $\Lambda^2_{L^2}=\Lambda^-_{L^2}\cong\Z$ and 
$\chi_{L^2}(M_V)=1$ and $\sigma_{L^2}(M_V)=-1$. Moreover for the 
hyper-K\"ahler metric $r=0$ and $s=0$ hence by (\ref{euler}) the exponent of 
${\rm Im}\:\tau /8\pi^2$ is
\[\frac{1}{240\pi^2}\int\limits_{M_V}\vert R\vert^2\dd V
=\frac{1}{30}\cdot\frac{1}{8\pi^2}
\int\limits_{M_V}{\rm tr}(R\wedge *R)=\frac{1}{30}\chi_{L^2}(M_V)=
\frac{1}{30}.\]
Consequently 
\[Z(M_V,g_V,\tau )
=\frac{(\det\:'\triangle_0)^2}{(\det\:'\triangle_1)^{\frac{1}{2}}}
\:\vartheta_{q^-_{M_V}}(-\overline{\tau})
\left(\frac{{\rm Im}\:\tau}{8\pi^2}\right)^{\frac{1}{30}}.\]
The modular weights (\ref{alfsulyok}) are 
$(-\frac{1}{30},\frac{7}{15})$. Since $H^1(W;\R )\cong H^1(S^3;\R )=\{ 0\}$ 
holds for the infinity of this space, condition (\ref{perem}) is 
satisfied therefore the modular anomaly can be removed by 
adding counter terms to the original naive action (\ref{hatas}).

Finally we remark that in our opinion the topological condition 
(\ref{perem}) is an artifact and should be removed from the construction 
by calculating $\zeta_{\triangle_k}(0)$ directly on the non-compact 
manifold without truncating it and imposing any boundary condition.

\end{document}